\documentclass[11pt,a4paper]{amsart}
\usepackage{amssymb,amsmath}
\usepackage{color}
\def\R{\mathbb{R}}
\def\m1{{I\!\!M}}
\newcommand{\rdue}{\R^2}

\newcommand{\elle}[1]{L^{#1}}

\newcommand{\grad}{\nabla}

\renewcommand{\to}{\rightarrow}
\newcommand{\pa}{\partial}

\renewcommand{\i}{\infty}

\newcommand{\dt}{\delta}

\newcommand{\al}{\alpha}

\newcommand{\om}{\Omega}

\newcommand{\lm}{\lambda}
\newcommand{\ov}[1]{\overline{#1}}
\newcommand{\uv}[1]{\underline{#1}}

\newcommand{\fo}{\forall}
\newcommand{\e}[1]{{\,\dsp e^{\dsp #1}}}
\newcommand{\Va}{\dfrac{V(x)}{|x|^{2\al}}}
\newcommand{\bo}{\mbox{\rm O}}

\newtheorem{theorem}{Theorem}[section]
\newtheorem{proposition}[theorem]{Proposition}
\newtheorem{lemma}[theorem]{Lemma}
\newtheorem{corollary}[theorem]{Corollary}
\newtheorem{remark}[theorem]{Remark}
\newtheorem{definition}[theorem]{Definition}
\newcommand{\finedim}{\hspace{\fill}$\square$}

\newcommand{\rife}[1]{(\ref{#1})}
\newcommand{\scp}{\scriptscriptstyle}
\newcommand{\sscp}{\scriptstyle}
\newcommand{\dsp}{\displaystyle}
\newcommand{\noi}{\noindent}
\renewcommand{\dfrac}{\displaystyle\frac}
\newcommand{\brm}{\begin{remark}\rm}
\newcommand{\erm}{\end{remark}}
\newcommand{\bdf}{\begin{definition}\rm}
\newcommand{\edf}{\end{definition}}
\newcommand{\bte}{\begin{theorem}}
\newcommand{\ete}{\end{theorem}}
\newcommand{\bpr}{\begin{proposition}}
\newcommand{\epr}{\end{proposition}}
\newcommand{\ble}{\begin{lemma}}
\newcommand{\ele}{\end{lemma}}
\newcommand{\bco}{\begin{corollary}}
\newcommand{\eco}{\end{corollary}}
\newcommand{\beq}{\begin{equation}}
\newcommand{\eeq}{\end{equation}}
\newcommand{\bdm}{\begin{displaymath}}
\newcommand{\edm}{\end{displaymath}}

\newcommand{\graf}[1]{\left\{\begin{array}{ll}#1\end{array}\right.}

\def\sideremark#1{\ifvmode\leavevmode\fi\vadjust{\vbox to0pt{\vss
 \hbox to 0pt{\hskip\hsize\hskip1em
 \vbox{\hsize2.1cm\tiny\raggedright\pretolerance10000
  \noindent #1\hfill}\hss}\vbox to15pt{\vfil}\vss}}}%

\begin{document}
\numberwithin{equation}{section}
\hfuzz=2pt
\frenchspacing

\title[Cosmic strings and Alexandrov's inequality]{Self gravitating cosmic strings and the Alexandrov's inequality for Liouville-type equations.}

\author[D. Bartolucci \& D. Castorina]{Daniele Bartolucci$^{(1,\ddag)}$, Daniele Castorina$^{(2)}$}

\thanks{2000 \textit{Mathematics Subject classification:} 35B45, 35J60, 35J99. }

\thanks{$^{(1)}$Daniele Bartolucci, Department of Mathematics, University
of Rome {\it "Tor Vergata"}, \\  Via della ricerca scientifica n.1, 00133 Roma,
Italy. e-mail: bartoluc@mat.uniroma2.it}

\thanks{$^{(2)}$Daniele Castorina, Dipartimento di Matematica, Universit\'a di Padova, \\
Via Trieste 63, 35121 Padova, Italy. e-mail: castorin@math.unipd.it}

\thanks{$^{(\ddag)}$Research partially supported by MIUR project {\sl
Metodi Variazionali ed Equazioni Differenziali non lineari.}}

\begin{abstract}
Motivated by the study of self gravitating cosmic strings, we
pursue the well known method by C. Bandle to obtain a weak version of the classical Alexandrov's isoperimetric inequality. In fact we
derive some quantitative estimates for weak subsolutions of a Liouville-type equation with conical singularities. Actually we succeed in generalizing 
previously known results, including Bol's inequality and pointwise estimates, to the case where the solutions solve the equation just
in the sense of distributions.  Next, we derive some \uv{new} pointwise estimates suitable to be applied to a class of singular cosmic 
string equations. Finally, interestingly enough,  we apply these results to establish a minimal mass property for solutions of the cosmic string equation 
which are \uv{supersolutions} of the singular Liouville-type equation.
\end{abstract}
\maketitle
{\bf Keywords}: Self gravitating cosmic strings, Isoperimetric inequalities on surfaces, Singular Liouville equations, Conical
singularities.
\\

\section{Introduction}
\noi Let $a>0$ and $N>-1$. Motivated by the study of cosmic strings configurations in the framework of 
Einstein's general relativity, the systematic study of the equation
\beq\label{1505.1}
-\Delta u = (e^{a u}+\left|x\right|^{2N} e^{u})=:f \quad\mbox{ in }\quad \R^2,\\
\eeq
has been initiated in some recent works \cite{cgs}, \cite{PoT}, \cite{Tar14}. A basic question about  \rife{1505.1} is for which values of 
$\beta:=\frac{1}{2\pi}\int_{\R^2} f$ it admits solutions. This is already a non trivial task for radial solutions, where however the sharp thresholds are known \cite{PoT}, while it is 
still open in the general case. A first step to set up this problem is to try to characterize the values of $\beta$ (which we denote by $\beta_0$) which can be achieved along 
blow up sequences. Let us assume for the moment that $x_0\in \R^2$ is a blow up point (see Definition \ref{blowuppt} below) 
relative to a sequence of solutions of \rife{1505.1}. By using
known arguments based on Brezis-Merle's estimates \cite{bm}, one can prove that $x_0$ is isolated and that 
$\beta_0\geq \min\{2(1+N_{-}),\frac2a\}$, where $N_{-}=\min\{0,N\}$, see \cite{Tar14}. However it is also well known that this "minimal mass" is not sharp in general, 
the optimal values being found via a refined analysis of the blow up behavior of a sequence of solutions in the same spirit of \cite{ls}, 
see \cite{Tar14}.  Concerning this point, the more subtle situation is met when the blow up point $x_0$ is not finite. In this case, after the Kelvin transform 
$x \mapsto \frac{x}{|x|}$ and even if $N>0$, one ends up with a \uv{singular} problem 
(see \rife{eqkelvin} below) where the weights multiplying the nonlinearities $e^{a u},\;e^{u}$ are power-type functions unbounded (in general) near the origin. To find out a 
minimal mass gets even more complicated in this case.
This is one of our motivations, since we succeed in obtaining a simpler and shorter argument yielding better (than $\min\{2(1+N_{-}),\frac2a\}$) 
estimates of the minimal masses which in some cases are also sharp, see Theorem \ref{infblowup}. 
Our result is based on some \uv{new} pointwise estimates suitable to be applied 
to \rife{1505.1} and \rife{eqkelvin}, see Proposition \ref{propexample}.\\

Moreover, as far as we are interested in the set of $\beta$ for which \rife{1505.1} admits 
solutions, then in general the above estimates on $\beta_0$ do not suffice to catch the optimal range. In this more general setting one needs to work out a much harder 
blow up analysis which takes into account the many allowed asymptotic behaviours as well as the sharp values of the corresponding local masses to be obtained via a local 
Pohozaev's identity. This kind of studies are based on the concentration-compactness-quantization theory for Liouville type equations \cite{bm}, \cite{yy}, \cite{ls} 
and its more recent generalizations to some singular cases, see \cite{bcct}, \cite{BM3}, \cite{bt}, \cite{T3} and \cite{T2} for a 
comprehensive exposition of the subject.  Interestingly enough, it turns out that, at least 
in a particular range of values for $a$, our result already provides the best possible estimate for $\beta$, i.e. $\frac4a$, see $(j)$ in Theorem \ref{infblowup}. 
Indeed, if $N>0$ and $\frac{1}{N+1}<a<\frac{2}{N+1}$, then the \underline{sharp} lower threshold 
(see \cite{Tar14}) reads $\beta \geq \frac4a$.\\

Next, let us fix some notations. Here and in the rest of this paper $\om\subset \rdue$ is any open, bounded and \uv{simply connected} domain, 
and $\widehat{V}\in L^\infty(\om) $ any measurable
function satisfying,
\beq\label{vi}
0<a\leq \widehat{V}\leq b<+\i,\;\mbox{for a.a. }\;x\in\om,
\eeq
for some fixed $0<a\leq b<\i$. We assume that $\al\in[0,1)$ and that $\om$ contains the origin, $0\in\om$, and set
$$
h_\al(x)=-2\al\log{|x|},\;\;x\in\R^2\setminus{\{0\}}.
$$
We are interested in the analysis of some quantitative properties of subsolutions of the singular
Liouville-type equation \cite{Lio}, \cite{Pic}
\beq\label{P}
-\Delta v = \Va\e{v}\quad\mbox{in}\quad\om,
\eeq
where
\beq\label{pesosub1}
V(x)=\widehat{V}(x) e^{g(x)},\quad g\in C^{0}(\om)\;
\mbox{and subharmonic in}\;\om.
\eeq
Unless otherwise specified we assume $\al\in(0,1)$ and with an abuse of notation also write $\frac1\al$
meaning $+\infty$ whenever $\al=0$.\\
We will be concerned with some integral inequalities and pointwise estimates for weak subsolutions of \rife{P}. 
However,  if $u$ is a solution of  \rife{1505.1}, then it is a \uv{supersolution} of \rife{P}. 
This is another interesting point about our result, since we are able to handle a case which does not fit in the general assumptions needed to apply 
the Alexandrov-Bol's inequality \cite{Band}.
We achieve this goal by a clever use of an auxiliary unknown which is in fact a subsolution of a singular Liouville equation, 
the inequality in \rife{P} being obtained by neglecting some negative terms. However, no informations are at hand concerning these negative terms other than 
their weak regularity properties, as dictated by the equation \rife{1505.1} itself. This is why, unlike previously known results, we have to keep weaker 
regularity assumptions about $v$ (see {\bf (a)} and {\bf (b)} below).\\
Actually, this is a general problem that arises in the study of solutions of singular cooperative Liouville systems \cite{PoT2}. 
In those cases one has to deal just with supersolutions of singular Liouville-type equations 
sharing poor regularity properties. Moreover, in that situation the above mentioned arguments based on the 
Pohozaev's identity fail to provide the needed sharp local masses. Another interesting point about our result is that it provides a direct method to obtain some 
estimates about the local masses in this case as well.\\ 

So, either one of the following assumptions will be made about $v$.\\
{\bf (a)} $v\in \elle{1}(\om)$,
$V e^{h_\al+v}\in \elle{1}(\om)$ and $v$ is a solution of \rife{P} in the sense of distributions.\\
By the results in \cite{bm} (see \cite{bm} Remarks 2 and 5) and standard elliptic regularity
theory \cite{Gilb} any such solution satisfies
\beq\label{rega}
v\in W_{\rm loc}^{2,p}(\om\setminus \{0\})
\cap W^{2,q}(\om),\;\mbox{for any}\;p\in[1,+\infty),\, q\in\left[1,\frac{1}{\al}\right).
\eeq
{\bf (b)} $v$ is a strong subsolution of \rife{P}, that is,
\beq\label{regb}
v\in W_{\rm loc}^{2,p}(\om\setminus \{0\}) \cap
W^{2,q}(\om),
\;\mbox{for some}\;p>2\;\mbox{and some}\;q\in\left(1,\frac{1}{\al}\right),
\eeq
and satisfies
\beq\label{subsol}
-\Delta v \leq \Va\e{v}\quad \mbox{for a.a. }x\in \; \om.
\eeq

Obviously, by the Sobolev embedding Theorem, we have $v\in C^{1}_{\rm loc}(\om\setminus\{0\}) \cap C^{0}(\om)$, a fact
that will be used throughout the discussion with no further comments.

\bdf\label{simp}{\it
We say that $\omega\subset \R^2$ is a \uv{simple} domain, if
it is an open and bounded domain whose boundary $\pa\omega$ is the support of
a rectifiable Jordan curve.
We will also say that $\omega\subset \R^2$ is a \uv{regular} domain if it is an
open and bounded domain whose boundary $\pa\omega$ is the union of finitely many rectifiable
Jordan curves. We will denote by $\ov{\omega_B}$ the closure of the (possibly
disconnected) bounded
component of $\R^2\setminus\omega$ and by $\omega_B$ its interior. The set of regular domains includes the set of simple domains which is
just characterized by the condition $\ov{\omega_B}=\emptyset$.}
\edf

Let $\omega\Subset \om$ be any regular domain and let us define,
$$
L_{\al}(\pa\omega)=\int\limits_{\pa\omega}{\dsp e}^{\frac{v+g+h_\al}{2}}d\ell,
$$
where $d\ell$ denotes the standard arc-length on $\pa\omega$,
\begin{equation}\label{malpha}
M_{\al}(\omega)=\int\limits_{\omega}{\dsp e}^{v+g+h_\al}dx,
\end{equation}

\beq\label{1003.1}
\widehat{K}(x)=\frac{1}{2}\widehat{V}(x),\;\;x\in\om,
\eeq
and, for any $K_0\geq 0$ and $\lm\in[0,1)$
\beq\label{1003.2}
\gamma_{\scp \omega}(\lm,K_0)=
2\pi\lm
\;\;-\!\!\!\!\!\int\limits_{\{\widehat{K}> K_0\}\cap\omega}\!\!\!(\widehat{K}-K_0){\dsp e}^{v+g+h_\al}dx.
\eeq
In \cite{Band} C. Bandle adopted a weighted rearrangement
argument to prove the following inequality,
\beq\label{120314.0.1}
L^2_{\al}(\pa\omega)\geq \left(2\gamma_{\scp \omega}(1-\al,K_0)-K_0M_{\al}(\omega)\right)M_{\al}(\omega),\quad
\mbox{with }\omega\mbox{ any simple domain},
\eeq
whenever $2\gamma_{\scp \omega}(1-\al,K_0)-K_0M_{\al}(\omega)>0$ and $0\in\omega$.
We remark that the inequality in \cite{Band} is more general since it allows one to replace $g+h_\al$ with the difference of two subharmonic 
(not necessarily continuous) functions, say $g_{+}-h_{-}$. In that case
$2\pi\al$ should be replaced by the total mass relative to the distributional Laplacian of $-h_{-}$ in $\omega$. 
Our proof could be modified to allow these kind of data. However, on one side the analysis of the cosmic strings equation \rife{1505.1} 
does not require this more general assumptions, which indeed fit more naturally in the framework of the geometric problem 
(see the short discussion below for more details about this point). On the other side, because of 
our weak formulation of the problem, the needed modifications would call up for a more involved technical discussion. Therefore we skip this part here for the sake of simplicity.   
Actually \rife{120314.0.1} is just a local formulation of a singular version of a classical isoperimetric inequality on surfaces, the Alexandrov's inequality,
which has a long history in geometry, \cite{Ale}, \cite{bol}, \cite{fia}, \cite{Hub0}. From this point of view, and in case
$g\equiv0$ and $\al=0$, then
$L_\al$ and $M_\al$ are just the local expressions of the length and the area of a portion of a surface whose
volume element (in local isothermal coordinates) takes the form $e^v dx$ and whose Gaussian curvature is
$\widehat{K}(x)$. If $\al\neq 0$ then the origin $\{0\}$ is a conical singularity of order $\al$, see \cite{Troy}, and the effect of a non vanishing and possibly 
singular $g$ is equivalent in general to the introduction of other kind of singularities, see \cite{Ale}, \cite{Band}. 
However, as mentioned above, our main motivation is not just with respect to the geometric problem, which is in fact
by now well understood, see \cite{Bur} (\S 2.1, \S 2.2), \cite{Band} (\S I.3.5)
and \cite{Oss}, \cite{Oss2} for more details and \cite{topp}, \cite{topp1} for more recent results
in this direction.
It is worth to remark that, besides the classical geometric applications \cite{barjga}, \cite{BDeM}, \cite{BdMM},  \cite{cama}, \cite{Troy}, other well known 
physical problems motivate this kind of studies, such as those arising in the statistical
mechanics description of guiding-centre plasmas \cite{BM3}, turbulent Euler flows \cite{BLin3}, \cite{clmp2}
and self-gravitating classical systems \cite{B2}, \cite{w}.\\ 

Actually, a  weaker (and well known in geometry) inequality, known as Bol's inequality \cite{bol}, which in our notations reads as
\beq\label{120314.0}
L^2_\al(\pa\omega)\geq\frac12  \left(8\pi(1-\al) - M_\al(\omega)\right)M_\al(\omega),\quad
\mbox{with }\omega\mbox{ any simple domain},
\eeq
newly derived in the P.D.E. setting in \cite{Band0}, has been since then widely used
in the analysis of elliptic problems with exponential nonlinearities in two-dimensions, see
\cite{Band}, \cite{B1}, \cite{barjga}, \cite{BM2}, \cite{cl}, \cite{ccl}, \cite{CLin3}, \cite{S}, \cite{suz}
and more recently \cite{bl}, \cite{BLin3}, \cite{BLT}.
The reason for its success in this context essentially relies on the fact that it can be used
to build up a weighted rearrangement-type argument \cite{Band}
to estimate the eigenvalues of linear Dirichlet problems on abstract surfaces,
or either of the linearized equations for problems with exponential nonlinearities such as \rife{P}.\\

One of our aims here is to prove (see Theorem \ref{Al} below) a weaker version of the Alexandrov's
inequality as derived in \cite{Band} \S I.3.5 where $v$ was assumed to be an analytic
subsolution of \rife{P}. Actually, analyticity in \cite{Band} was a reasonable and well suited assumption as it allowed the author to handle 
the general framework described above. On the contrary, as far as we are concerned with \rife{pesosub1}, it seems a rather strong
assumption. This will be our first concern since we are not aware
of any proof of \rife{120314.0.1} under weaker assumptions such as those in {\bf (a)} or {\bf (b)}.
Actually, in \cite{suz}, \rife{120314.0} is proved for $V \equiv 1$, $\al=0$ and $v\in C^{2}(\om)\cap C^{0}(\,\ov{\om}\,)$
a classical subsolution of \rife{P}.
Other results in \cite{ccl} and more recently in \cite{bl}, \cite{BLin3}, where $v$ is still assumed to be of class $v\in C^{2}(\om)\cap C^{0}(\,\ov{\om}\,)$, 
are concerned with various forms of 
Bol's inequality. So we also obtain a unified proof and a generalization of these inequalities in case $V$ satisfies \rife{pesosub1} and $\al\in(0,1)$, see Theorem \ref{Al-bdy} and
Corollary \ref{Boli}.\\

We also make a point which deserves a separate discussion.
The case where $\omega\Subset \om$ is multiply connected in \rife{120314.0} with $V\equiv1$ and $\al=0$,
has been handled in Lemma 4.2 in \cite{ccl} by using the fact that
the inequality holds on simply connected domains, and by assuming that
\beq\label{120314.1}
\int\limits_{\om} e^v\leq 8\pi.
\eeq
However, it is readily seen that the same argument will not
work when trying to extend the Alexandrov's inequality \rife{120314.0.1}
on multiply connected domains $\omega\Subset \om$ from simply connected ones. The point is that,
even in case $g\equiv0$ and $\al=0$, the Alexandrov's inequality is much stronger than Bol's inequality,
since it provides a kind of measure of "how far" the isoperimetric ratio is from the one with $\widehat{K}\equiv K_0$ assumed to be constant.
Therefore one of our goals is to show that there is no need to assume \rife{120314.1},
because in fact not only \rife{120314.0} but also the Alexandrov's inequality hold on multiply connected domains $\omega\subset\om$ as well,
see Theorem \ref{Al} and Corollary \ref{Boli} below. It is understood that
this observation applies as far as $\om$ itself is assumed to be simply connected, otherwise these inequalities
are well known to be false in general, see \cite{Bur}, \cite{Oss2} and more recently \cite{BLin3}. Although we will not discuss it here, 
it is worth to remark that a major improvement in the P.D.E. analysis of the Bol's inequality
in case \uv{$\om$ is assumed to be multiply connected} has been
recently obtained in \cite{BLin3}.\\

We divide the proof of Theorem \ref{Al} into two steps. In the first step we manage to
pass from the elliptic inequality \rife{subsol} to an elliptic equation for an auxiliary function.
In the second one we work out a weighted rearrangement argument. The advantage of this approach is that the arguments in these two steps
are suitable to be moved as they stand in the proof of the needed pointwise estimates, see Theorem \ref{ptwise} below. Even in this case
we obtain a generalization of other results \cite{Band}, \cite{suz} previously derived under stronger assumptions.
In any case these pointwise estimates will be the main tool in the analysis of the blow up phenomenon relative to \rife{1505.1}.\\
We remark that part of our results are likely to be well known to experts. For example a
weaker form of Corollary \ref{Boli} has been already used in \cite{B1} but its proof was given for granted there.\\

This paper is organized as follows.
In section \ref{s0} we discuss an inequality by A. Huber and its consequences.
In section \ref{s1} we prove the Alexandrov's inequality (Theorem \ref{Al}).
Some Corollaries and improved versions of Theorem \ref{Al} are discussed in section \ref{s2}.
The pointwise estimates are discussed in section \ref{pointwise}. Finally the
application arising in cosmic strings theory is discussed in section \ref{example}. A technical result is established in the Appendix.\\

\textit{Acknowledgement:} We would like to thank Prof. G. Tarantello for many helpful discussions 
and in particular for bringing to our attention the application of those pointwise estimates to the analysis of the cosmic strings equation.

\section{Preliminaries: An inequality by A. Huber and its consequences.}\label{s0}

In this section we discuss a particular form of the Huber's inequality \cite{Hub} suitable to be applied to our problem.
The result due to A. Huber \cite{Hub} is more general than the one we discuss here and for the sake of simplicity we report 
only the part which will be needed throughout.\\
Let $\omega \subset \R^2$ be a regular domain and let us keep the same notations as in the introduction above.
Then we have
\bte\label{hub-teo}[Huber's inequality, see \cite{Hub}]  Let $\omega\Subset \R^2$ be a simple domain. Then it holds:
\beq\label{Hubineq-sing}
\left(\,\int\limits_{\pa \omega}{\dsp e}^{\frac{g+h_\al}{2}}d\ell\right)^2
\geq 4\pi\left(1 -\al\right)\int\limits_{\omega}{\dsp e}^{g+h_\al}dx,\,\,\,\mbox{if}\;\;0\in\ov{\omega},
\eeq
\beq\label{Hubineq-reg}
\left(\,\int\limits_{\pa \omega}{\dsp e}^{\frac{g+h_\al}{2}}d\ell\right)^2
\geq 4\pi\int\limits_{\omega}{\dsp e}^{g+h_\al}dx,
\,\,\mbox{if}\;\;0\notin\ov{\omega}.
\eeq

\ete

We will need the following generalization of Huber's result.

\bte\label{hub-teo-2} Let $\omega\Subset \R^2$ be a regular domain. Then it holds:\\
\beq\label{Hubineq-sing-1}
\left(\,\int\limits_{\pa \omega}{\dsp e}^{\frac{g+h_\al}{2}}d\ell\right)^2
\geq 4\pi\left(1 -\al\right)\int\limits_{\omega}{\dsp e}^{g+h_\al}dx,
\,\,\mbox{if}\;\;0\in \ov{\omega} \cup \ov{\omega_B},
\eeq
\beq\label{Hubineq-reg-1}
\left(\,\int\limits_{\pa \omega}{\dsp e}^{\frac{g+h_\al}{2}}d\ell\right)^2
\geq 4\pi\int\limits_{\omega}{\dsp e}^{g+h_\al}dx,
\,\,\mbox{if}\;\;0\notin \ov{\omega} \cup \ov{\omega_B},
\eeq
where $d\ell$ denotes the arc-length on $\pa \omega$.\\
In particular, if $\omega$ is not simply connected, then all the inequalities are strict.

\ete
\brm{\it
In view of this result we will be free
to use the inequality \rife{Hubineq-sing-1} in all the cases considered so far.}
\erm

\proof
In view of Theorem \ref{hub-teo} we are obviously left to discuss the cases where $\omega$ is not simply connected
and prove in particular that in all those cases the inequalities are strict.\\
Let us assume for the moment that $\omega=\omega_1\setminus\ov{\omega_0}$ for a pair of
simple domains such that $\omega_0\subset\subset \omega_1$ and
$\pa\omega = \pa\omega_1 \cup \pa \omega_0$.
For any domain $\omega\subset\rdue$,
let us set
$$
\ell(\pa\omega)=\;\int\limits_{\pa \omega} {\dsp e}^{\frac{g+h_\al}{2}}d\ell,\qquad
m(\omega)=\int\limits_{\omega} {\dsp e}^{g+h_\al}dx.
$$
Thus, if $0\in\ov{\omega}$, we may use \rife{Hubineq-sing} and \rife{Hubineq-reg} to obtain
$$
\ell^2(\pa\omega)=\ell^2(\pa\omega_1 \cup \pa \omega_0)>
\ell^2(\pa\omega_1)+\ell^2(\pa\omega_0)\geq
$$
$$
4\pi(1-\al) m(\omega_1)+ 4\pi (1-\al) m(\omega_0)>4\pi(1-\al)  m(\omega),
$$
and the desired conclusion follows in this case. On the other side, if $0\notin\ov{\omega_1}$,
we obtain
$$
\ell^2(\pa\omega)=\ell^2(\pa\omega_1 \cup \pa \omega_0)>\ell^2(\pa\omega_1)+\ell^2(\pa\omega_0)\geq
$$
$$
4\pi m(\omega_1)+ 4\pi m(\omega_0)>
4\pi m(\omega).
$$
Finally, if $0\in\ov{\omega_0}$ we obtain
$$
\ell^2(\pa\omega)=\ell^2(\pa\omega_1 \cup \pa \omega_0)>\ell^2(\pa\omega_1)+\ell^2(\pa\omega_0)\geq
$$
$$
4\pi(1-\al) m(\omega_1)+ 4\pi(1-\al) m(\omega_0)>
4\pi(1-\al) m(\omega).
$$
The case where $\R^2\setminus\omega$ has finitely many bounded components readily follows by an
induction argument on the number of "holes" of $\omega$.\finedim

\bigskip

\section{Alexandrov's inequality}\label{s1}
In this section we prove a weak version of the celebrated Alexandrov's inequality
in the form first provided in \cite{Band} for analytic subsolutions of \rife{P}.
We remark that, as in \cite{Band}, the case where the equality holds
in Alexandrov's inequality can be completely characterized as well. We keep the same notations as in Definition
\ref{simp}.

\bte\label{Al}[Alexandrov's inequality \cite{Ale}, \cite{Band}, \cite{Bur}, \cite{Oss2}]\\
Let $\om\subset \rdue$ be any open, bounded and simply connected domain with $0\in\om$ and
$\omega\Subset\om$ be any relatively compact regular domain.
Fix $\al\in[0,1)$ and let $v$ satisfy either {\bf (a)} or {\bf (b)}. Then:\\
if $0\in{\ov{\omega}\cup \ov{\omega_B}}$, then, for any $K_0\geq 0$, it holds
\beq\label{Alex}
L^2_{\al}(\pa\omega)\geq \left(2\gamma_{\scp \omega}(1-\al,K_0)-K_0M_{\al}(\omega)\right)M_{\al}(\omega),
\eeq
while
if $0\notin\ov{\omega}\cup \ov{\omega_B}$, then, for any $K_0\geq 0$,
it holds
\beq\label{Alex-2}
L^2_{\al}(\pa\omega)\geq \left(2\gamma_{\scp \omega}(1,K_0)-K_0M_{\al}(\omega)\right)M_{\al}(\omega).
\eeq
\ete

\brm{\it One would be tempted to include the strict inequality
$K_0 M_{\al}(\omega)<2\gamma_{\scp \omega}(1-\al,K_0)$
as an hypothesis to obtain \rife{Alex}, at least since otherwise \rife{Alex} is trivially satisfied.
However that inequality is not needed during the proof. The same consideration holds for
\rife{Alex-2} and the inequality $K_0 M_{\al}(\omega)<2\gamma_{\scp \omega}(1,K_0)$. In particular there is even no
need to assume the nonnegativity of $\gamma_{\scp \omega}(\lm,K_0)$. However it is well known that
if $K_0 M_{\al}(\omega)\geq 2\gamma_{\scp \omega}(1-\al,K_0)$, then in general no bound for $M_{\al}(\omega)$ is possible
in terms of $L^2_{\al}(\pa\omega)$, see for example \cite{Oss} p.1207.}
\erm

\brm 
{\it If $\omega$ is not simply connected and $0\in \omega_B$, then there is no hope to obtain \rife{Alex-2}, that is, 
\rife{Alex} with $\al=0$. In other words, the fact that the a "hole" contains the singularity makes the isoperimetric ratio 
lower by an amount that is at least as larger as the full weight of the singularity. Actually this is the 
same reason why Alexandrov's inequality fails on a general multiply connected domain $\om$, see \cite{BLin3} p.14 for further details.}
\erm

\brm\label{remLap}{\it We recall some well known facts and refer the reader to \cite{Gilb} \S 2.8 for further details. 
An open domain $\omega$ is said to be \uv{regular with respect to the Laplacian} if
for any $x\in\pa\omega$ there exist a barrier function at $x$ relative to $\omega$. In particular any regular 
(according to Definition \ref{simp}) domain $\omega$ is regular with respect to the Laplacian. 
We will use the fact that the classical Dirichlet problem with continuous boundary data on a 
regular domain admits a unique solution in the class of functions that are continuous up to the boundary.}
\erm

\noindent {\bf Proof.}\\
Once \rife{Alex} is proved for $0\in \omega\cup\omega_B$, a straightforward approximation
argument shows that indeed \rife{Alex} holds whenever $0\in \ov{\omega}\cup\ov{\omega_B}$ as well.
Also, once the result has been established for $K_0>0$, then the case $K_0=0$ is worked out by
an elementary limiting argument. Therefore we will just discuss the cases $0\in \omega\cup\omega_B$ and
$K_0>0$.\\
The proofs of either \rife{Alex-2} or of \rife{Alex} in case $\al=0$ as well as that of
\rife{Alex} in case $0\in \omega_B$ are easier than that of \rife{Alex} in case $0\in \omega$ and $\al\in(0,1)$.
Therefore, in particular to avoid repetitions, we will just be concerned
with the proof of \rife{Alex} in case $0\in \omega$ and $\al\in(0,1)$.\\
In view of \rife{rega}, it is obvious that if  $v$ satisfies {\bf (b)} then it also satisfies {\bf (a)}, so we will
be concerned just with the former case.\\

\noi
{\bf The Proof of \rife{Alex} when $v$ satisfies (b) and $0\in\omega$, $\al\in(0,1)$ and $K_0>0$.}\\
Unless otherwise specified, we assume $0\in\omega$, $\al\in(0,1)$ and $K_0>0$.\\
{\bf Step 1.}\\
For any fixed relatively compact and regular domain $\omega\Subset\om$
such that $0\in\omega$ we can find
an open, simply connected and smooth domain $\om_0$ such that
$$
\omega \Subset \om_0\Subset \om.
$$
Let us also define
$$
f(x):=-\Delta v - \Va\e{v},\quad  x\in \om\setminus\{0\},
$$
which in view of \rife{subsol} and {\bf (b)} satisfies
$$
f\in L^{p}_{\rm loc}(\om\setminus\{0\})\cap L^{q}_{\rm loc}(\om),
\; \text{ and } f(x)\leq 0, \quad\mbox{for a. a. } x\in \om,
$$
for some $p>2$ and some $q\in\left(1,\frac{1}{\al}\right)$.
Therefore, in view of \rife{regb}, Theorem 9.15, Corollary 9.18 and Lemma 9.17 in \cite{Gilb} we see that the linear problem
\beq\label{0403.10}
-\Delta w = f(x)\;\; \mbox{in}\;\;\om_0,
\qquad w = 0 \;\; \mbox{on}\;\; \pa \om_0,
\eeq
admits a unique solution $w \in W^{2,p}_{\rm loc}(\om_0\setminus\{0\})\cap W^{2,q}(\om_0)\cap C^{0}(\,\ov{\om_0}\,)$,
for some $p>2$ and some $q\in\left(1,\frac{1}{\al}\right)$. Clearly $w$ is subharmonic  (see \cite{Gilb} \S 2.8 and Ex. 2.7, 2.8).\\
Next let $g_1$ be the Perron's (see \cite{Gilb} \S 2.8) solution of
$\Delta g_1=0$ in $\omega$, $g_1=v$ on $\pa\omega$. Since $v\in C^{0}(\overline{\omega})$, then
in view of Remark \ref{remLap} $g_1$ is well
defined and continuous up to the boundary. Finally let us set $\eta=v-w-g_1$.
Then, since $v$ satisfies \rife{regb}, we see that
$\eta\in  W^{2,p}_{\rm loc}(\omega\setminus\{0\})\cap W^{2,q}(\omega)\cap C^{0}_{0}(\ov{\omega})$
for some $p>2$ and some $q\in\left(1,\frac{1}{\al}\right)$ and satisfies
\beq\label{eqeta-sub}
-\Delta \eta= 2\widehat{K} \e{\sscp g_2+h_\al}\e{\eta}\;\; \mbox{for a. a. }x\in\; \omega,
\quad \eta=0\;\; \mbox{on}\;\; \pa\omega,
\eeq
with
\beq\label{subh-sub}
g_2=w+g_1+g\;\mbox{continuous in}\;\ov{\omega}\; \mbox{and subharmonic in}\; \omega.
\eeq
Thus, in particular by using \rife{regb}, we conclude that
\beq\label{reg.1b-sub}
\eta\in W^{2,p}_{\rm loc}(\omega\setminus\{0\})\cap C^{1}_{\rm loc}(\omega\setminus\{0\})\cap
W^{2,q}(\omega)\cap C^0_0(\ov{\omega}),
\eeq
for some $p>2$ and some $q\in\left(1,\frac{1}{\al}\right)$.
Since $\eta \in W^{2,q}(\omega)$, by Sobolev embedding Theorem, standard elliptic estimates
(see Lemma 9.17 in \cite{Gilb}), an elementary smoothing argument and the maximum principle for $W^{2,2}(\omega)$ solutions
(see Theorem 9.1 in \cite{Gilb}), we see that $\eta\geq 0$ in $\omega$. In particular,
since $\eta \in W^{2,2}_{\rm loc}(\omega\setminus\{0\})$,
a well known version of the strong maximum principle (see Theorem 9.6 in \cite{Gilb}) applies
and we see that $\eta$ is also strictly positive in $\omega\setminus\{0\}$. Finally,
due to the fact that $\eta$ is also weakly superharmonic, then it satisfies (see \cite{Gilb} Ex. 2.8)
$$
\eta(x)\geq \frac{1}{2\pi r}\int\limits_{\pa B_r(x)}\eta \, d\ell,
$$
for any $x\in \omega$ whenever $B_r(x)\Subset \om$. Whence we also have $\eta(0)>0$ and we conclude that
\beq\label{2.71bis}
\eta(x)>0\quad \fo\, x\in \omega\quad\mbox{and}\quad \eta(x)=0\;\iff\;x\in\pa\omega.
\eeq

\bigskip

\noi
{\bf Step 2.}
Setting $t_m=\max\limits_{\ov{\omega}} \eta$ and
$$
d\tau={\dsp e}^{g_2+h_\al}dx,\;\;d\sigma={\dsp e}^{\frac{g_2+h_\al}{2}}d\ell,
$$
for any $t\in [0, t_m)$ we define
\beq\label{defG}
\om(t)=\{x\in\omega\,|\,\eta(x)>t\},\quad \Gamma(t)=\{x\in\omega\,|\,\eta(x)=t\},\quad
\mu(t)=\int\limits_{\om(t)}d\tau,
\eeq
being understood that $\Gamma$ is defined for $t=t_m$ as well.\\
Since $\eta$ satisfies \rife{eqeta-sub} then $\Gamma(t)$ has null two-dimensional measure,  whence  we
conclude that $\mu$ is continuous. Moreover, in view of \rife{2.71bis}, the following relations hold true
\beq\label{reg.2-sub}
\om(0)=\omega,\qquad \Gamma(0)=\pa\omega,\qquad \mu(0)=\int\limits_{\omega}d\tau.
\eeq
The fact that $\Gamma(t)$ has null 2-dimensional measure is well known.\\
Clearly we can extend $\mu$ on $[0,t_m]$ by setting $\mu(t_m)=\lim\limits_{t\nearrow t_m}\mu(t)=0^{+}$
whence $\mu\in C^{0}([0,t_m])$. Next, since $\eta$ satisfies \rife{eqeta-sub}, then
by a well known consequence (see for example \cite{bz} p.158) of the co-area
formula and of Sard's Lemma we see that

\beq\label{070214.1-sub}
\frac{d\mu(t)}{dt}=-\int\limits_{\Gamma(t)}\frac{{\dsp e}^{\sscp g_2+h_\al}}{|\nabla \eta|}\,d\ell,\;\;
\mbox{for a. a. }t\in [0,t_m].
\eeq

\brm\label{remclv}{\it
Let $I\subseteq [t_0,t_m]$ be the set where \rife{070214.1-sub} holds.
Obviously we may assume that if $t\in I$ then $t$ is not a critical level of $\eta$. In particular,
letting $t_0\in(0,t_m]$ be the unique $t$ such that $0\in \Gamma(t_0)$, there is no loss of generality
in assuming $t_0\notin I$.}
\erm

For any $s\in[0,\mu(0))\equiv [\mu(t_m),\mu(0))$, we introduce the following weighted
decreasing rearrangement of $\eta$,
\beq\label{etadef-sub}
\eta^*(s)= |\{t\in [0,t_m]:\mu(t)>s\}|,
\eeq
where $|E|$ denotes the Lebesgue measure of a Borel set $E\subset \R$. Setting
$\eta^*(\mu(0))=0$, and using the fact that $\Gamma(t)$ has null 2-dimensional measure, it is not
difficult to check that $\eta^*\in C^{0}([0,\mu(0)])$ is the inverse of $\mu$ on $[0,t_m]$ (hence continuous in $[0,\mu(0)]$),
coincides with the distribution function of $\mu$ and it is also strictly decreasing and
differentiable almost everywhere. A crucial point at this stage is to prove
that $\eta^*$ is not just continuous but also locally absolutely continuous.
It turns out that in fact it is even locally Lipschitz in $(0,\mu(0))$, see Appendix \ref{App}. In particular, in view of
\rife{070214.1-sub}, we obtain
\beq\label{070214.1-sub-1}
\frac{d\eta^*(s)}{ds}=-
\left(\,\int\limits_{\Gamma(\eta^*(s))}\frac{{\dsp e}^{\sscp g_2+h_\al}}{|\nabla \eta|}\,d\ell\right)^{-1},
\eeq
for any $s\in I^{*}$, where $[0,\mu(0)]\setminus I^{*}$ is a set of null measure and
$\eta^{*}(I^{*})=I$, $\mu(I)=I^*$.\\
Next, let us define
$$
F(s)=2K_0\int\limits_{\om(\eta^{*}(s))}{\dsp e}^{\eta}d\tau,\quad s\in[0,\mu(0)],
$$
where
\beq\label{fin-sub1}
F(\mu(0))=2K_0\int\limits_{\omega}{\dsp e}^{\eta}d\tau=2K_0M_\al(\omega),
\eeq
and we have set
\beq\label{fin-sub0}
F(0)=\lim\limits_{s\searrow 0^+}F(s)=0.
\eeq
Clearly $F(s)$ is strictly increasing, continuous, and even locally Lipschitz in $(0,\mu(0))$.
In fact it satisfies
$$
|F(s)-F(s_0)|\leq C |\mu(\eta^{*}(s))-\mu(\eta^{*}(s_0))|= C|s-s_0|,\quad \fo\,0=\mu(t_m)<s_0<s<\mu(0),
$$
for a suitable constant $C>0$. In particular it holds
$$
\int\limits_{\om(\eta^*(s))}{\dsp e}^{\dsp u}d\tau=
\int\limits_{0}^{s}{\dsp e}^{\dsp\eta^*(\lm)}d\lm,\;\;\;\fo\,s\in[0,\mu(0)],
$$
so that
\beq\label{F2}
\frac{dF(s)}{ds}= 2K_0e^{\eta^*(s)},\quad
\frac{d^2 F(s)}{ds^2}=2K_0\frac{d\eta^*(s)}{ds}\,{\dsp e}^{\dsp \eta^*(s)}=\frac{d\eta^*(s)}{ds}\frac{d F(s)}{ds},\quad\fo\,s\in I^*.
\eeq
For any $s\in I^*$ the Cauchy-Schwartz inequality yields,
\beq\label{281109.1.1}
\left(\,\int\limits_{\Gamma(\eta^*(s))} d\sigma\right)^2\leq
\left(\,\int\limits_{\Gamma(\eta^*(s))}\frac{{\dsp e}^{g_2+h_\al}}{|\nabla \eta|}\,d\ell\right)
\left(\,\int\limits_{\Gamma(\eta^*(s))}|\nabla \eta|d\ell\right)=
\eeq
$$
\left(-\frac{d\eta^*(s)}{ds}\right)^{-1}
\left(\,\int\limits_{\Gamma(\eta^*(s))}\left(-\frac{\pa \eta}{\pa \nu}\right)d\ell\right),
$$
where $\nu=\frac{\nabla \eta}{|\nabla \eta|}$ is the exterior unit normal to $\om(\eta^*(s))$ and we have used
\rife{070214.1-sub-1}. Since $\eta$ satisfies \rife{reg.1b-sub}, then it is easy to check that

$$
\int\limits_{\Gamma(\eta^*(s))}\left(-\frac{\pa \eta}{\pa \nu}\right)d\ell=
\int\limits_{\om(\eta^*(s))}2K \e{\eta}d\tau,
$$
for any $s\in I^*$, whence, in particular we deduce that
$$
\int\limits_{\Gamma(\eta^*(s)))}\left(-\frac{\pa \eta}{\pa \nu}\right)d\ell\leq
2\int\limits_{\{K> K_0\}\cap\om(\eta^*(s)))}(K-K_0)\e{\eta}d\tau+
2K_0\int\limits_{\om(\eta^*(s)))}\e{\eta}d\tau\leq
$$
$$
2\int\limits_{\{K> K_0\}\cap\omega}(K-K_0)\e{\eta}d\tau+
2K_0\int\limits_{\om(\eta^*(s)))}\e{\eta}d\tau=
2\tilde{\gamma}^+_\omega(K_0) + F(s),
$$
for any $s\in I^*$, where
\begin{equation}\label{gammatilde}
\tilde{\gamma}^+_\omega(K_0)=\int\limits_{\{K> K_0\}\cap\omega}(K-K_0)\e{\eta}d\tau.
\end{equation}
Plugging this estimate in \rife{281109.1.1} we conclude that,
\beq\label{iso-1-sub}
\left(\,\int\limits_{\Gamma(\eta^*(s))} d\sigma\right)^2
\leq \left(-\frac{d\eta^*(s)}{ds}\right)^{-1}\left[2\tilde{\gamma}^+_\omega(K_0)+F(s)\right],
\eeq
for any $s\in I^*$. Since $\omega$ is open and by assumption $0\in\omega$, by setting $t_0=\eta(0)$ we have that
either $t_0<t_m$, and then $0\in \om(t)$ for any $t\in(0,t_0)$
and $0\notin \om(t)$ for any $t\in[t_0,t_m)$,
or $0\in \om(t)$ for any $t\in(0,t_m)$. The discussion concerning the latter case is easier
so, in order to avoid repetitions, we will only prove \rife{Alex} for $t_0<t_m$.
Let $s_0\in (0,\mu(0))$ satisfy $\eta^*(s_0)=t_0$.
Clearly we can apply Huber's inequality \rife{Hubineq-sing-1} to conclude that,
\beq\label{Huber-sub}
\left(\,\int\limits_{\Gamma(\eta^*(s))} d\sigma\right)^2\geq
4\pi(1-\al)\mu(\eta^*(s))=4\pi(1-\al)s,\,\forall\;s\in I^*\cap (s_0,\mu(0)).
\eeq
On the other side, using the maximum principle and the fact that $\Gamma(t)$ has null two-dimensional measure,
it is not difficult to check that $\R^2\setminus \om(t)$ has no bounded components with $t\in(t_0,t_m)$. Therefore
the singular point $0\in\R^2$ cannot be contained in a bounded component of
$\R^2\setminus \om(t)$ with $t\in(t_0,t_m)$, so we can use \rife{Hubineq-reg-1}
to conclude that
\beq\label{Huber-2-sub}
\left(\,\int\limits_{\Gamma(\eta^*(s))} d\sigma\right)^2\geq 4\pi \mu(\eta^*(s))\geq 4\pi(1-\al)\mu(\eta^*(s))=4\pi(1-\al)s,
\,\forall\;s\in I^*\cap (0,s_0).
\eeq

\noi
Substituting the last two inequalities in \rife{iso-1-sub}, we obtain,
$$
4\pi(1-\al)s\leq\left(-\frac{d\eta^*(s)}{ds}\right)^{-1}\left[2\tilde{\gamma}^+_\omega(K_0)+F(s)\right]
,\,\forall\;s\in I^*,
$$
Multiplying by $\frac{dF(s)}{ds}\left(-\frac{d\eta^*(s)}{ds}\right)$
 it is readily seen that the last two inequalities are equivalent to,
$$
4\pi(1-\al)s \frac{dF(s)}{ds}\left(\frac{d\eta^*(s)}{ds}\right)
+2\frac{dF(s)}{ds}\tilde{\gamma}^+_\omega(K_0)+\frac{dF(s)}{ds}F(s)\geq 0,\,
\forall\;s\in I^*,
$$
and we conclude that
\begin{equation}\label{pjalfa} 
\frac{d}{ds}\left[4\pi(1-\al)s \frac{dF(s)}{ds}-4\pi(1-\al)F(s)+
2\tilde{\gamma}^+_\omega(K_0)F(s)+
\frac{1}{2}(F(s))^2 \right]\geq 0,\,\forall\;s\in I^*.
\end{equation}

Let $P(s)$ denote the function within square brackets in \eqref{pjalfa}.
Since $F$ and $\eta^*$ are both continuous and
locally Lipschitz continuous in $[0,\mu(0)]$ and, in view of \rife{F2},
since $\frac{dF(s)}{ds}$ is continuous and
locally Lipschitz continuous in $[0,\mu(0)]$ too,
then we come up with the inequality
$$
P(\mu(0))-P(0)\geq 0.
$$
Therefore we can use \rife{fin-sub1}, \rife{fin-sub0}, \rife{F2} and $\eta^*(\mu(0))=0$ to obtain
$$
4\pi(1-\al)\mu(0)+
(2\tilde{\gamma}^+_\omega(K_0)-4\pi(1-\al))M_\al(\omega)+K_0 M^2_\al(\omega)\geq 0.
$$
By the Huber's inequality \rife{Huber-sub} once more and \rife{reg.2-sub} we have
\beq\label{24-basta-sub}
L^2_{\al}(\pa\omega)=\left(\,\int\limits_{\pa\omega}{\dsp e}^{\frac{v+\widehat{g}+h_\al}{2}}d\ell\right)^2=
\left(\,\int\limits_{\Gamma(0)} d\sigma \right)^2\geq 4\pi(1-\al)\mu(0)\geq
\eeq
$$
(4\pi(1-\al)-2\tilde{\gamma}^+_\omega(K_0))M_\al(\omega)-K_0 M^2_\al(\omega)=
$$
$$
(2(2\pi(1-\al)-\tilde{\gamma}^+_\omega(K_0))-K_0 M_\al(\omega))M_\al(\omega)=
(2\gamma_{\scp \omega}(1-\al,K_0)-K_0M_{\al}(\omega))M_{\al}(\omega),
$$
which concludes the proof of \rife{Alex}.
\finedim

\bigskip
\section{Corollaries and/or improved versions of Theorem \ref{Al}}\label{s2}

The first slightly improved version of Theorem \ref{Al} we wish to discuss has to do with the
possibility to allow $\omega=\om$. Actually this was the point of view in \cite{suz} where however
$\om$ was assumed to be smooth, $v\in C^{2}(\om)\cap C^{0}(\ov{\om})$ and $\al=0$ and $V\equiv 1$.
Here we have,

\bte\label{Al-bdy}
Let $\om\subset \rdue$ be a simple domain of class $C^{1,1}$ (see \cite{Gilb})
such that $0\in\om$ and fix $\al\in[0,1)$.  Let $g$ (as defined in \rife{pesosub1}) satisfy $g\in C^0(\ov{\om})$ and for any $r>0$ small enough $v$ satisfy
\rife{subsol} and
\beq\label{regb-bdy}
v\in W^{2,p}(\om\setminus B_r(0))\cap
W^{2,q}(\om)\cap C^0(\ov{\om}),
\;\mbox{for some}\;p>2\;\mbox{and for some}\;q\in\left(1,\frac{1}{\al}\right).
\eeq
Then \rife{Alex} holds with $\omega\equiv \om$.
\ete
\proof
The proof follows step by step the one of Theorem \ref{Al} with very few and minor
modifications, so we just skip it.
\finedim

\bigskip

Next we discuss a Corollary which is an oversimplified
(we choose $V\equiv 1$) \uv{and weaker} (we fix $K_0=\frac12$) version of Theorem \ref{Alex}.
It improves and generalizes former results in \cite{Band0},
\cite{ccl} and in \cite{suz}.
When $\al=0$ it is just the Bol's inequality, first
discovered by G. Bol \cite{bol}.

\bco\label{Boli}[Singular Bol's inequality \cite{bol} and \cite{Band}, \cite{ccl}, \cite{suz}]
Let $\om\subset \rdue$ be any open, bounded and simply connected domain with $0\in\om$ and
$\omega\Subset\om$ be any relatively compact regular domain. Fix $\al\in[0,1)$ and
let $v$ satisfy either {\bf (a)} or {\bf (b)} with $\widehat{V}\equiv 1$. Then:\\
{\bf (i)} if $0\in\ov{\omega}\cup\ov{\omega_B}$, then it holds
\beq\label{Alex-suz-1003}
L^2_{\al}(\pa\omega)\geq\frac12  \left(8\pi(1-\al) - M_{\al}(\omega)\right)M_{\al}(\omega),
\eeq
while if $0\notin\ov{\omega}\cup\ov{\omega_B}$, then it holds
\beq\label{Alex-2-suz-1003}
L^2_{\al}(\pa\omega)\geq \frac12 \left(8\pi -  M_{\al}(\omega)\right)M_{\al}(\omega);
\eeq

{\bf (ii)} if $\om$ is assumed to be simple and of class $C^{1,1}$,  $\widehat{V}\equiv 1$, and if
 $v$ satisfies \rife{regb-bdy} and \rife{subsol}, then \rife{Alex-suz-1003} holds with
$\omega\equiv \om$.
\eco
\proof
In view of \rife{vi}, \rife{1003.1} and \rife{1003.2}, if $\widehat{V}\equiv 1$ and $K_0=\frac12$, then
$\widehat{K}\equiv \frac12$ and we have
$$
\gamma_{\scp \omega}\left(1-\al,\frac12\right)=2\pi(1-\al)
\;\;-\!\!\!\!\!\int\limits_{\{\widehat{K}>\frac12\}\cap\omega}\!\!\! \left(\widehat{K}-\frac12\right){\dsp e}^{v+g+h_\al}dx\equiv
2\pi(1-\al),
$$
so that
$$
2\gamma_{\scp \omega}(1-\al,K_0)-K_0M_{\al}(\omega)=4\pi(1-\al)-\frac12M_{\al}(\omega).
$$
Therefore \rife{Alex-suz-1003}, \rife{Alex-2-suz-1003} 
readily follow from \rife{Alex}, \rife{Alex-2}.\\
The inequality claimed in {\bf (ii)} follows by the same argument and Theorem \ref{Al-bdy}.
\finedim

\bigskip

\section{Pointwise estimates based on the weighted rearrangement}\label{pointwise}

In this section we will establish pointwise estimates for strong subsolutions. More specifically we will prove the following theorem, which improves and generalizes former 
results in \cite{Band0}, \cite{suz}:

\begin{theorem}\label{ptwise}
Let $\om\subset \rdue$ be any open, bounded and simply connected domain with $0\in\om$ and
$\omega\Subset\om$ be any relatively compact regular domain. Fix $\al\in[0,1)$ and
let $v$ satisfy either {\bf (a)} or {\bf (b)} with $\widehat{V}\equiv 1$. Then:\\
{\bf (i)}
If $0\in\ov{\omega}\cup\ov{\omega_B}$ and $M_\al(\omega) < 8\pi(1-\al)$ then it holds
\beq\label{eq-ptwise1}
\max_{\overline{\omega}} e^{v} \leq \left( 1 - \frac{M_\al(\omega)}{8\pi(1-\al)} \right)^{-2} \max_{\partial \omega} e^{v},
\eeq
while if $0\notin\ov{\omega}\cup\ov{\omega_B}$ and $M_\al(\omega) < 8\pi$, then it holds
\beq\label{eq-ptwise2}
\max_{\overline{\omega}} e^{v} \leq \left( 1 - \frac{M_\al(\omega)}{8\pi} \right)^{-2} \max_{\partial \omega} e^{v}.
\eeq
{\bf (ii)} If $\om$ is assumed to be simple and of class $C^{1,1}$ and if
 $v$ satisfies \rife{regb-bdy} and \rife{subsol}, then \rife{eq-ptwise1} holds with
$\omega\equiv \om$.
\end{theorem}

\proof
We will be mainly concerned with \eqref{eq-ptwise1}, which is the more involved case. 
The rest of the proof goes along the same lines with minor changes, and we skip this part here to avoid repetitions.

Let $\eta$ be exactly the same function as defined in Step 1 in the proof of Theorem \ref{Al}. We recall from Step 2 that if we set $\widehat{V} \equiv 1$ and $K_0 = \frac12$, then
$\widehat{K} \equiv \frac12$, $F(s)=\int_{\om(\eta^{*}(s))} e^{\eta}\, d\tau$, $s\in[0,\mu(0)]$ and $\tilde{\gamma}^+_\omega (1/2) = 0$ (see \eqref{gammatilde}). Moreover, setting
\begin{equation}\label{eq-ptwise3}
P_{\al} (s) = 4\pi s F'(s)-4\pi F(s)+\frac{1}{2}(F(s))^2, \, \forall \; s \in (0,\mu(0)) 
\end{equation}
then from \eqref{pjalfa} we have 
\begin{equation}\label{eq-ptwise3bis}
P_{\al} (s) \geq P_{\al} (0) = 0  \; \forall \, s \in [0,\mu(0)]. 
\end{equation}
We also recall that if we define
\begin{equation}\label{eq-ptwise4}
J_\al(s) = \frac{s}{F(s)} - \frac{s}{8\pi(1-\al)}, \quad  s \in (0,\mu(0))
\end{equation}
then we can differentiate and from \eqref{eq-ptwise3}, \eqref{eq-ptwise3bis} obtain  
\begin{equation}\label{eq-ptwise5}
J_\al '(s) = -\frac{P_{\al} (s)}{4\pi(1-\al)F^2(s)} \leq 0, \quad s \in (0,\mu(0))
\end{equation}
Hence $J_\al$ is $C^1$ and nonincreasing. Thanks to l'Hopital rule we can extend it by continuity to $s = 0$ as follows
\begin{equation}\label{eq-ptwise6}
\lim_{s \to 0^+} J_\al(s)= \lim_{s \to 0^+} \frac{s}{F(s)} = \lim_{s \to 0^+} \frac{1}{F'(s)} = \frac{1}{F'(0)}.
\end{equation}
In particular, since $F'(0) = \max_{\overline{\omega}} e^\eta$, for any $s \in (0,\mu(0))$ we get
\begin{equation}\label{eq-ptwise7}
\max_{\overline{\omega}} e^\eta \leq \frac{1}{J_\al(s)}.
\end{equation}
On the other hand, again by \eqref{eq-ptwise3bis}, for $s \in (0,\mu(0))$ we deduce
\begin{equation}\label{eq-ptwise8}
s F'(s) \geq F^2(s) \left( \frac{1}{F(s)} - \frac{1}{8\pi(1-\al)} \right).
\end{equation}
Notice that thanks to the hypothesis \eqref{eq-ptwise1} we have $F(s) \leq F(\mu(0)) = M_\al (\omega) < 8\pi(1-\al)$,
which in turn by \eqref{eq-ptwise8} and $F' >0$ implies that for $s \in (0,\mu(0))$ it holds
\begin{equation}\label{eq-ptwise9}
\begin{split}
&J_\al(s) = s \left( \frac{1}{F(s)} - \frac{1}{8\pi(1-\al)} \right) = \frac{sF'(s)}{F'(s)} \left( \frac{1}{F(s)} - \frac{1}{8\pi(1-\al)} \right)\\
&\geq \frac{F^2(s)}{F'(s)} \left( \frac{1}{F(s)} - \frac{1}{8\pi(1-\al)} \right)^2 =
\frac{1}{F'(s)} \left( 1 - \frac{F(s)}{8\pi(1-\al)} \right)^2.
\end{split}
\end{equation}
Thus combining \eqref{eq-ptwise7} and \eqref{eq-ptwise9}, we obtain that for any $s \in (0,\mu(0))$ it holds
\begin{equation}\label{eq-ptwise10}
\max_{\overline{\omega}} e^\eta \leq F'(s) \left( 1 - \frac{F(s)}{8\pi(1-\al)} \right)^{-2}.
\end{equation}
So, passing to the limit as $s \nearrow \mu(0)$ in \eqref{eq-ptwise10} and noticing that $F'(\mu(0))= e^{\eta (\mu(0))} = 1$ and $F(\mu(0)) = M_\al (\omega)$, we finally have
\begin{equation}\label{eq-ptwise11}
\max_{\overline{\omega}} e^\eta \leq \left( 1 - \frac{M_\al (\omega)}{8\pi(1-\al)} \right)^{-2}.
\end{equation}
Now, since $\eta = v-g_3$ (with $g_3$ subharmonic and continuos in $\omega$, see Step 1) and since $\eta = 0$ on $\partial \omega$,
then by the weak maximum principle and \eqref{eq-ptwise10} we obtain
\begin{equation}\label{eq-ptwise12}
\begin{split}
&\max_{\overline{\omega}} e^v = \max_{\overline{\omega}} e^{\eta+g_3} \leq \max_{\overline{\omega}} e^\eta \max_{\overline{\omega}} e^{g_3}
\leq \left( 1 - \frac{M_\al (\omega)}{8\pi(1-\al)} \right)^{-2} \max_{\overline{\omega}} e^{g_3}\\
&= \left( 1 - \frac{M_\al (\omega)}{8\pi(1-\al)} \right)^{-2} \max_{\partial \omega} e^{g_3} =
\left( 1 - \frac{M_\al (\omega)}{8\pi(1-\al)} \right)^{-2} \max_{\partial \omega} e^v,
\end{split}
\end{equation}
which is \eqref{eq-ptwise1} as claimed.
\finedim

\bigskip
\bigskip

\section{An application of Alexandrov-Bol's inequality arising in cosmic strings theory.}\label{example}

The aim of this section is to discuss an explicit example arising in cosmic strings theory, see \cite{Tar14}. 
To simplify the exposition, in this section we let $N_-=\min\{0,N\}$, while $B_\dt$ will be used to denote a ball of radius $\dt>0$ centered at an arbitrarily fixed point $x_0\in\R^2$. 
We have the following:
\begin{proposition}\label{propexample}
For any $N,L >-1$ and $a>0$, let $u \in L^1_{\rm loc}(\R^2)$ be a solution (in the sense of distributions) of:
\begin{equation}\label{eqexample0tris}
\begin{cases}
-\Delta u = (|x|^{2N} e^{u}+ |x|^{2L} e^{a u}):= f \quad\mbox{ in }\quad \R^2,\\
\int_{\R^2} f < +\infty.
\end{cases}
\end{equation}
If $M_{a,N,L} := \max\{1,a\} \int_{B_\delta} f < 8\pi(1+\min\{N_-, L_-\})$, then
\begin{equation}\label{eqexample1tris}
\max_{\overline{B}_\delta} e^{u} \leq \left( 1 - \frac{M_{a,N,L}}{8\pi(1+\min\{N_-, L_-\})} \right)^{-2} \max_{\partial B_\delta} e^{u}.
\end{equation}
\end{proposition}

\brm {\it
The proof of this result relies on Theorem \ref{ptwise} and the fact that suitably modified logarithms of the datum $f$ satisfy ${\bf (b)}$. We are indebted with G. Tarantello 
for bringing this fact to our attention.}    
\erm

\begin{proof}
We will discuss the case $N\geq L$ first, and then indicate how to derive the latter case $N < L$ from the former.\\
\noi
{\bf The Proof of \rife{eqexample1tris} in case $N\geq L$.}\\
If $N \geq L$ we rewrite \eqref{eqexample0tris} as follows:
$$-\Delta u = |x|^{2L} (e^{a u}+ |x|^{2(N-L)} e^{u}):= |x|^{2L} g \quad\mbox{ in }\quad \R^2.$$
Observe that by the results in \cite{bm} we know that $u^+ \in L^{\infty} (B_{2\dt})$, which spells that $f \in L^\infty (B_{2\dt})$ whenever
$L\geq0$, while $f \in L^\infty(\om_1)$ on any compact subset $\om_1\subset \ov{ B_{2\dt}}\setminus \{0\}$ whenever $-1 < L < 0$. 
By standard elliptic regularity $u$ is then smooth away from the origin if $-1 < L < 0$, while it is of class $C^{2}(B_{2\dt})$ if $L\geq 0$. 
Thus we shall further divide our discussion in two subcases, namely $L \geq 0$ and $-1<L<0$.\\
\noi
{\bf The Proof of \rife{eqexample1tris} in case $N\geq L$ and $L\geq0$.}\\
If $L\geq 0$ the estimate \eqref{eqexample1tris} is a direct consequence of \eqref{eq-ptwise1} with $\al=0$.
In fact we will prove that the function
\begin{equation}\label{eqexample2}
\eta_a =  \log(g) + \log(\max\{1,a\}),\quad x\neq 0,
\end{equation}
satisfies ${\bf (b)}$ with $V =|x|^{2L}$, and $\al=0$, that is,
\beq\label{040414.1}
\eta_a \in W_{\rm loc}^{2,p}(B_{2\dt} \setminus \{0\}) \cap W^{2,q}(B_{2\dt})
\eeq
for some $p>2$ and some $q>1$ and
\begin{equation}\label{eqexample3}
-\Delta \eta_a \leq |x|^{2L} e^{\eta_a} \quad \text{ a.e. in } \R^2.
\end{equation}
Clearly, as far as we are concerned with \rife{040414.1}, and since here $L \geq 0$, then
it will be enough to prove the $W^{2,q}(B_{2\dt})$ regularity of $\eta_a$ for some $q>1$.
At this point, let us set
\begin{equation}\label{eqexample4}
V_{a}(x)=\left|x\right|^{2(N-L)}e^{(1-a)u}
\end{equation}
and
\begin{equation}\label{eqexample5}
\xi_{a}(x):=u (x)+\dfrac{1}{a}\log\left(1+V_{a}(x)\right),
\end{equation}
so that
\begin{equation}\label{eqexample6}
-\Delta u= |x|^{2L} \left(1+V_{a}(x)\right)e^{au}= |x|^{2L} e^{a\xi_a}.
\end{equation}
In order to compute $\Delta\xi_a$, we observe that, since $N\geq L$, then $1+V_a(x)$ is bounded away from zero. Hence we
have
\begin{equation}\label{eqexample7}
\Delta\log\left(1+V_{a}(x)\right)=\dfrac{\Delta V_{a}(x)}{1+V_{a}(x)}-
\dfrac{\left|\nabla V_{a}(x)\right|^{2}}{\left(1+V_{a}(x)\right)^{2}},
\end{equation}
and, for $x\neq 0$,
\begin{equation}\label{eqexample8}
\Delta V_{a}=(1-a)V_{a}(x)\Delta u+\dfrac{\left|\nabla V_{a}(x)\right|^{2}}{V_{a}(x)}.
\end{equation}
So, by using  \eqref{eqexample6}, \eqref{eqexample7} and \eqref{eqexample8}, for $x\neq 0$, we find
\begin{equation}\label{eqexample9}
\Delta\log\left(1+V_{a}(x)\right)=\left(\dfrac{(1-a)V_{a}(x)}{1+V_{a}(x)}\right)\Delta u+
\dfrac{\left|\nabla V_{a}(x)\right|^{2}}{V_{a}(x)\left(1+V_{a}(x)\right)^{2}}.
\end{equation}
Notice that the r.h.s. in \eqref{eqexample9} may be singular at $x=0$, but  since $u \in C^2 (\overline{B_{2\dt}})$,
then the asymptotic behavior is dictated just by the term $|x|^{2(N-L)}$ in the definition of $V_a (x)$ as given by \eqref{eqexample4}.
In particular, since $V_{a}(x) = \bo(|x|^{2(N-L)})$ and $|\nabla V_{a}(x)| = \bo(|x|^{2(N-L)-1})$ as $x \approx 0$,
we conclude that:
\begin{equation}\label{eqexample9bis}
\dfrac{\left|\nabla V_{a}(x)\right|^{2}}{V_{a}(x)\left(1+V_{a}(x)\right)^{2}} =
\graf{\bo \left(|x|^{2(N-L)-2} \right),\quad N>L \\ \bo \left(1 \right),\quad N=L}.
\end{equation}
Thus, from \eqref{eqexample9bis}, we see that the r.h.s. of \eqref{eqexample9} is in $L^q (B_{2\dt})$ for some $q>1$ whenever
$N\geq L$,  and the required regularity of $\eta_a$ is established. Next, in view of \eqref{eqexample5} and \eqref{eqexample9} we obtain,
\begin{equation}\label{eqexample10}
\begin{aligned}
-\Delta\xi_a&=-\Delta u\left(\dfrac{a+V_{a}(x)}{a\left(1+V_{a}(x)\right)}\right)-\dfrac{\left|\nabla V_{a}(x)\right|^{2}}{aV_{a}(x)\left(1+V_{a}(x)\right)^{2}}\\
&\leq |x|^{2L} e^{a\xi_a}\left(\dfrac{a+V_{a}(x)}{a\left(1+V_{a}(x)\right)}\right),\quad\mbox{ a.e. in }\mathbb{R}^{2}.
\end{aligned}
\end{equation}
Finally we distinguish two cases: if $a\in(0,1]$, from \eqref{eqexample10} we find,

$$-\Delta\xi_a\leq  \dfrac{|x|^{2L}}{a}\left(\dfrac{a+V_{a}(x)}{1+V_{a}(x)}\right)e^{a\xi_a}\leq  \dfrac{|x|^{2L}}{a}e^{a\xi_a}.$$

In particular, letting $\eta_{a}:= a\xi_a$ then $\eta_a$ satisfies \eqref{eqexample3}. On the other hand, for $a>1$,
the previous computations imply that,

$$-\Delta\xi_a\leq |x|^{2L} \left(\dfrac{1+\dfrac{1}{a}V_{a}(x)}{1+V_{a}(x)}\right)e^{a\xi_a}\leq |x|^{2L} e^{a\xi_a},$$
and the desired conclusion again follows for $\eta_a:=a\xi_a+\log a$.\\
\noi
{\bf The proof of \rife{eqexample1tris} in case $N \geq L$ and $-1<L<0$.}\\
The proof for $-1<L<0$ is inspired by the previous but is slightly more delicate. This time the estimate \eqref{eqexample1tris} will be a direct consequence of \eqref{eq-ptwise1} with $\al=-L$. In fact we will prove that in this case $\eta_a$, defined as in \eqref{eqexample3}, satisfies ${\bf (b)}$ with $V \equiv C|x|^{2L}$, and $\al=-L$, that is,
$$
\eta_a \in W_{\rm loc}^{2,p}(B_{2\dt} \setminus \{0\}) \cap W^{2,q}(B_{2\dt})
$$
for some $p>2$ and some $q>1$ and
\begin{equation}\label{eqexample3s}
-\Delta \eta_a \leq \frac{1}{|x|^{2|L|}}e^{\eta_a} \quad \text{ a.e. in } \R^2.
\end{equation}
Once again, it will be enough to prove the $W^{2,q}(B_{2\dt})$ regularity of $\eta_a$ for some $q>1$.
With the same definitions of $V_{a}$ and $\xi_{a}$ as given by \eqref{eqexample4} and \eqref{eqexample5}, 
notice that, since $u^+ \in L^{\infty} (B_{4\dt})$ and $N \geq L$, then $g$ is bounded.
The difference now is that the r.h.s of \eqref{eqexample6}, due to the presence of $|x|^{2L}$ with $L < 0$, is not bounded but belongs to $L^{q} (B_{4\dt})$ for any $1<q<\frac{1}{|L|}$,
which in turn yields that $u \in W^{2,q} (B_{2\dt})$. As before, by using \eqref{eqexample7}, for $x\neq 0$, we find
\begin{equation}\label{eqexample9s}
\Delta\log\left(1+V_{a}(x)\right)=\left(\dfrac{(a-1)V_{a}(x)}{1+V_{a}(x)}\right)\Delta u+
\dfrac{\left|\nabla V_{a}(x)\right|^{2}}{V_{a}(x)\left(1+V_{a}(x)\right)^{2}}.
\end{equation}
Notice that the r.h.s. in \eqref{eqexample9s} is singular at $x=0$. However, the first term can be estimated as follows:
$$\left|\dfrac{(a-1)V_{a}(x)}{1+V_{a}(x)} \Delta u\right| \leq |1-a| |\Delta u|$$
with $\Delta u \in L^{q} (B_{2\dt})$ for any $1<q<\frac{1}{|L|}$, while the second term can be treated exactly as in \eqref{eqexample9bis}. 
Thus, we see that the r.h.s. of \eqref{eqexample9s} is in $L^q (B_{2\dt})$ for some $q>1$ whenever $-1<L<0$, and the required regularity of $\eta_a$ is established. Next, in view of
\eqref{eqexample5} and \eqref{eqexample9s} we obtain,
\begin{equation}\label{eqexample10s}
\begin{aligned}
-\Delta\xi_a&=-\Delta u\left(\dfrac{1+a V_{a}(x)}{1+V_{a}(x)}\right)-\dfrac{\left|\nabla V_{a}(x)\right|^{2}}{V_{a}(x)\left(1+V_{a}(x)\right)^{2}}\\
&\leq \frac{1}{|x|^{2|L|}} e^{\xi_a}\left(\dfrac{1+aV_{a}(x)}{1+V_{a}(x)}\right),\quad\mbox{ a.e. in }\mathbb{R}^{2}.
\end{aligned}
\end{equation}
Once again we distinguish two cases: if $a\in(0,1]$, from \eqref{eqexample10s} we find,

$$-\Delta \xi_a \leq \frac{1}{|x|^{2|L|}} \left(\dfrac{1+a V_{a}(x)}{1+V_{a}(x)}\right) e^{\xi_a}\leq \frac{1}{|x|^{2|L|}} e^{\xi_a}.$$
In particular $\eta_{a}:= \xi_a$ satisfies \eqref{eqexample3s}. On the other hand, for $a>1$,
the previous computations imply,

$$-\Delta\xi_a\leq \frac{a}{|x|^{2|L|}} \left(\dfrac{\dfrac{1}{a} + V_{a}(x)}{1+V_{a}(x)}\right) e^{\xi_a}\leq \frac{a}{|x|^{2|L|}} e^{\xi_a},$$
and the desired conclusion again follows for $\eta_a:=\xi_a+\log a$.\\

Next, let us discuss the case $N < L$.\\
\noi
{\bf The proof of \rife{eqexample1tris} in case $N < L$.}\\
By a direct calculation we see that if $u$ satisfies \eqref{eqexample0tris} with parameters $(N,L,a)$, than the function $v=au$ satifies:
$$-\Delta v = a(|x|^{2N} e^{\frac{1}{a} v}+ |x|^{2L} e^{v}) = a|x|^{2N}(e^{\frac{1}{a} v} + |x|^{2(L-N)} e^{v}) \quad\mbox{ in } \quad \R^2$$
which is essentially \eqref{eqexample0tris} but with the different choice of parameters $(L,N,1/a)$. Furthermore we notice that:
$$M_{1/a,L,N} = \max\left\{1,\frac{1}{a}\right\} \int_{B_\delta} a(|x|^{2L} e^{v}+ |x|^{2N} e^{\frac{1}{a} v}) =$$
$$=\max\{1,a\} \int_{B_\delta} (|x|^{2L} e^{v}+ |x|^{2N} e^{\frac{1}{a} v}) = \max\{1,a\} \int_{B_\delta} (|x|^{2L} e^{au}+ |x|^{2N} e^{u}) = M_{a,N,L}$$
In particular, since $L>N$, this means that we can repeat the same discussion of the previous case with minor changes. Indeed, it will enough to replace $u$ and $L$ with $v$ and $N$ respectively, but with the choice of $V(x) = a |x|^{2L}$ (as in \eqref{eq-ptwise1}) and $V_a (x) = \left|x\right|^{2(L-N)}e^{(1 - \frac{1}{a})v}$ (as in \eqref{eqexample4}).
\end{proof}

\bigskip

\noi
At this point we can apply our results to the blow up analysis of self-gravitating strings \cite{Tar14}. For any $N>-1$ and $a>0$ we let $u \in L^1_{\rm loc}(\R^2)$ be a solution (in the sense of distributions) of:
\begin{equation}\label{eqexample0}
\begin{cases}
-\Delta u = (e^{a u}+\left|x\right|^{2N} e^{u}):= f \quad\mbox{ in }\quad \R^2,\\
\int\limits_{\R^2} f < +\infty.
\end{cases}
\end{equation}
The above problem \eqref{eqexample0} fits in the framework of Proposition \ref{propexample} for $L=0$, hence if $M_{a,N} := \max\{1,a\} \int\limits_{B_\dt} f<8\pi(1+N_-)$, we have that
\begin{equation}\label{eqexample1}
\max_{\overline{B_\delta}} e^{u} \leq \left( 1 - \frac{M_{a,N}}{8\pi (1+N_-)} \right)^{-2} \max_{\partial B_\delta} e^{u}.
\end{equation}
As a direct consequence of \eqref{eqexample1} one can find a lower threshold for the mass, as we shortly describe below, which is usually obtained 
through a longer and more complicated blow-up analysis, see \cite{Tar14}. Our approach allows for a shorter proof based on Proposition \ref{propexample}.\\

However we should take into account also the presence of blow-up at \textit{infinity}. In this regard, an important tool is the so called {\it Kelvin transform} $\hat{u}$ of $u$, given explicitly by:
\begin{equation}\label{kelvin}
\hat{u} (x) := u \left( \frac{x}{|x|^2} \right) + \beta_a \log \frac{1}{|x|}
\end{equation}
with $\beta_a = \frac{1}{2\pi} \int_{\mathbb{R}^2} e^{a u}+\left|x\right|^{2N} e^{u}$. It is not difficult to prove (see \cite{Tar14}) that if $u$ satisfies \eqref{eqexample0}, then $\hat{u}$ satisfies
\begin{equation}\label{eqkelvin}
-\Delta \hat{u} = |x|^{a \beta_a -4} e^{a \hat{u}} + |x|^{\beta_a - 2(N+2)} e^{\hat{u}} \quad \text{ in }\quad  \mathbb{R}^{2}
\end{equation}
where $\beta_a > \max\{2/a,2(N+1) \}$. In this way we see that equation \eqref{eqkelvin} is a particular case of \eqref{eqexample0tris} with $N=\frac{\beta - 2(N+2)}{2}$ and $L=\frac{a \beta-4}{2}$, which explains the full generality we pursued in Proposition \ref{propexample}. 

In order to specify our results, let us then consider a sequence $\left\{u_k\right\} \subset L^1 (\mathbb{R}^2)$ satisfying \eqref{eqexample0tris},
\begin{equation}\label{uk}
\begin{cases}
-\Delta u_k=e^{au_k}+\left|x\right|^{2N_k}e^{u_k}=:f_k\\
\beta_{a,k}:=\dfrac{1}{2\pi}\displaystyle{\int_{\mathbb{R}^{2}}e^{au_k}+\left|x\right|^{2N_k}e^{u_k}}
\end{cases}
\end{equation}
with
\begin{equation}\label{betauk}
N_k \to N > -1 \quad\text{ and }\quad \beta_{a,k} \to\beta\quad\mbox{ as }\quad k\to+\infty.
\end{equation}
Furthermore, let $\hat{u}_k$ be the sequence of corresponding Kelvin transforms as in \eqref{kelvin}, which satisfy
\begin{equation}\label{hatuk}
-\Delta \hat{u}_k = |x|^{a \beta_{a,k} -4} e^{a \hat{u}_k} + |x|^{\beta_{a,k} - 2(N_k +2)} e^{\hat{u}_k}:= \hat{f}_k \quad \text{ in } \mathbb{R}^{2},
\end{equation}
We give the following

\begin{definition}\label{blowuppt}
We say that $x_0 \in \mathbb{R}^2$ is a \textbf{blow-up point} for a sequence of functions $u_k$ if there exists a sequence 
$x_n \to x_0$ and a subsequence $k_n$ such that $u_{k_n} (x_n) \to +\infty$ as $n \to \infty$. 
We say that $u_k$ has a \textbf{blow-up point at infinity} (and we write $x_0 = \infty$) if the sequence $\hat{u}_{k}$ has a blow-up point at $x_0 = 0$.
\end{definition}

In the spirit of \cite{bm} and thanks to Proposition \ref{propexample} we will now see that a blow--up point, 
finite or not, requires a fixed amount of the $L^{1}$--norm of $f_k$.  In fact,
\begin{theorem}\label{infblowup}
Let $\left\{u_k\right\}_{k \in \mathbb{N}}$ be a sequence that satisfies \eqref{uk} and \eqref{betauk}. 
Suppose that $u_k$ has a blow-up point at $x_0$ and let $u_{k_n}$ be a subsequence as in Definition \ref{blowuppt}. The following alternatives hold:\\

$(j)$ If $x_0 \ne 0, \infty$, then for any $\delta>0$ it holds

\begin{equation}\label{almenofinito}
\beta_{\infty}:=\displaystyle{\liminf_{n\rightarrow+\infty}}\frac{1}{2\pi}\int_{B_{\delta}} f_{k_n} \geq \frac{4}{\max\left\{1,a\right\}}.
\end{equation}

$(jj)$ If $x_0 = 0$, then for any $\delta>0$ it holds

\begin{equation}\label{almenofinito0}
\beta_{\infty}:=\displaystyle{\liminf_{n\rightarrow+\infty}}\frac{1}{2\pi}\int_{B_{\delta}} f_{k_n} \geq \frac{4 (1+N_{-})}{\max\left\{1,a\right\}}.
\end{equation}
 
$(jjj)$ If $x_0 = \infty$, then for any $\delta>0$ it holds

\begin{equation}\label{almeno}
\beta_{\infty} := \displaystyle{\liminf_{n\rightarrow+\infty}}\frac{1}{2\pi}\int_{B_{\delta}} \hat{f}_{k_n} \geq 
\dfrac{(2+\min\{(a \beta -4)_{-}, (\beta - 2(N+2))_{-}\})}{\max\left\{1,a\right\}}.
\end{equation}

\end{theorem}

\brm{\it 
As mentioned in the introduction, the minimal mass value obtained in $(a)$ is sharp, see \cite{Tar14}.}
\erm

\begin{proof}
The proof of $(j)$ and $(jj)$ is an immediate consequence of \eqref{eqexample1}, 
since it can be shown (see \cite{Tar14}) that blow-up points are isolated. In fact, for each blow-up point $x_0$ one can find $\delta_0 = \delta_0 (x_0)$ 
such that for any $0 < \delta_1 < \delta_0$ it holds:
$$\sup_{|x|=\delta_1} u_{k_n} \leq C_{\delta_1}, \quad \forall \; n \in \mathbb{N}.$$
Hence, if for some $\delta>0$ either \eqref{almenofinito} or \eqref{almenofinito0} were not true, 
then, since $f_k$ is nonnegative, \eqref{eqexample1} would contradict the fact that $x_0$ is a blow-up point. 
The same argument applies to the proof of $(jjj)$. In this case we just use \eqref{eqexample1tris} of Proposition \ref{propexample} and \rife{betauk}, \rife{hatuk}. 
\end{proof}

\bigskip
\bigskip

\section{Appendix}\label{App}

In this appendix we prove that $\eta^{*}$ (as defined in \rife{etadef-sub}) is
locally Lipschitz in $(0,\mu(0))$.

\ble\label{lem-liploc}
For any $\al\in [0,1)$ and $0<\overline{a}\leq a<b\leq \overline{b}<\mu(0)$, there exist
$\ov{C}=\ov{C}(\overline{a},\overline{b},\al)>0$ such that
\beq\label{liploc}
\eta^*(a)-\eta^*(b)\leq \ov{C} (b-a).
\eeq
\ele
\proof
For fixed $a<b$ as above we can find an open set $\om_{a,b}$ such that
$$
\{x\in\omega\,:\, \eta^{*}(b)\leq \eta(x)\leq \eta^{*}(a) \}\Subset \om_{a,b}\Subset\omega.
$$
Using Green's representation formula, and in view of \rife{reg.1b-sub},
it is not difficult to check that
$$
|\grad \eta(x)|\leq C +C\int\limits_{\om_{a,b}}\frac{dy}{\dsp |x-y||y|^{2\al}},\;\forall \;x\in \ov{\om_{a,b}}.
$$
It is well known that
$$
\int\limits_{\om_{a,b}}\frac{dy}{\dsp |x-y||y|^{2\al}}\leq \int\limits_{B_R(0)}\frac{dy}{\dsp |x-y||y|^{2\al}}
\leq C(1+|x|^{1-2\al}),
$$
for some $R\geq 1$ depending on $\om_{a,b}$.
So we find
\beq\label{AppB-1}
|\grad \eta(x)|\leq C+C |x|^{1-2\al}\;\forall \;x\in \ov{\om_{a,b}},
\eeq
for a suitable constant $C>0$. Let 
$$
d\tau={\dsp e}^{g_2+h_\al}dx,
$$
with $g_2$ as in \rife{subh-sub}. Thus, by the co-area formula and Sard's Lemma we obtain
$$
b-a=\mu(\eta^*(b))-\mu(\eta^*(a))=
\int\limits_{\eta>\eta^*(b)} d\tau-\int\limits_{\eta>\eta^*(a)} d\tau=
\int\limits_{\eta^*(b)<\eta\leq\eta^*(a)} \hspace{-.5cm}d\tau\geq
$$
$$
\int\limits_{\eta^*(b)<\eta<\eta^*(a)} \hspace{-.5cm}d\tau=
\int\limits_{\eta^*(b)}^{\eta^*(a)}
\left(\, \int\limits_{\Gamma(t)}\frac{{\dsp e}^{g_2+h_\al}}{|\nabla \eta|}\, d\ell\right)dt\geq
\int\limits_{\eta^*(b)}^{\eta^*(a)}
\left(\, \int\limits_{\Gamma(t)}
\frac{{\dsp e}^{g_2+h_\al}}{C+C |x|^{1-2\al}}\, d\ell\right)dt\geq
$$
$$
C_1\int\limits_{\eta^*(b)}^{\eta^*(a)}
\left(\, \int\limits_{\Gamma(t)}
\frac{{\dsp e}^{g_2}}{|x|^{2\al} +|x|}\, d\ell\right)dt\geq
C_2\int\limits_{\eta^*(b)}^{\eta^*(a)}
\left(\, \int\limits_{\Gamma(t)}d\ell\right)dt=C_2\int\limits_{\eta^*(b)}^{\eta^*(a)}
\mathcal{L}_1(\Gamma(t))dt\geq
$$
$$
C_2\inf\limits_{\eta^{*}(b)\leq t\leq\eta^{*}(a) }\mathcal{L}_1(\Gamma(t))
\int\limits_{\eta^*(b)}^{\eta^*(a)}dt=\ov{C}(\eta^*(a)-\eta^*(b)),
$$
where $\mathcal{L}_1(\Gamma)$ is the 1-dimensional Lebesgue measure of the set $\Gamma$ and
in the last inequality we have used the standard isoperimetric inequality to conclude that
$$
\mathcal{L}_1(\Gamma(t))\geq 4\pi |\,\om(t)|\geq 4\pi |\,\om(\eta^*(\,\overline{a}\,))|>0,
$$
for any $\eta^{*}(b)\leq t\leq\eta^{*}(a)$.
\finedim

\end{document}